\documentclass[12pt,twoside]{article}
\usepackage[centertags]{amsmath}
\usepackage{amsfonts}
\usepackage{amssymb}
\usepackage{amsthm}
\usepackage{newlfont}
\usepackage{color}

\date{\color{green}  2013 November}
\def \n {\noindent}
\usepackage{fancyhdr}
\newcommand{\hf}{\hfill $\diamondsuit$}
\begin{document}

\begin{center}
{\Large{\color{blue}{\bf On regularized trace formula of Gribov semigroup genrated by the Hamiltonian of reggeon field theory in Bargmann representation}}}\\
\end{center}

\begin{center}
{\bf Abdelkader INTISSAR  $^{(*)}$ $^{(**)}$}\\
\end{center}

\begin{center}
\scriptsize{(*)Equiped'Analyse spectrale, UMR-CNRS n: 6134, Universit\'e de
Corse, Quartier Grossetti, 20 250 Cort\'e-France \\
 T\'el: 00 33 (0) 4 95 45 00 33\\
 Fax: 00 33 (0) 4 95 45 00 33\\
 e.mail:intissar@univ-corse.fr\\
\parskip = 3pt
(**)Le Prador,129 rue du commandant Rolland, 13008 Marseille-France}
\end{center}

\quad\\

\begin{center}

\fbox{\rule[-0.4cm]{0cm}{1cm} \n {\color{red}{\bf Abstract}}}\\

\end{center}

{\it \small{In J. Math. Anal. App. 305. (2005), we have considered the Gribov operator\\ $H_{\lambda'} = \lambda' S + H_{\mu,\lambda}$ acting on Bargmann space where
$ S = a^{*2} a^{2}$ and $ H_{\mu,\lambda} =\mu a^* a + \\i \lambda a^* (a+a^*)a$ with $i^{2} = -1$.\\
Here $a$ and $a^{*}$ are the standard Bose annihilation and creation operators satisfying the commutation relation $[a, a^{*}] = I$. In Reggeon field theory, the real parameters $\lambda{'}$ is the four coupling of Pomeron, $\mu$  is Pomeron intercept, $\lambda$ is the triple coupling of Pomeron and $i^{2} = -1$.\\
\quad We have given an approximation of the semigroup $e^{-tH_{\lambda'}}$ generated by the operator $H_{\lambda'}$. In particulary, we have obtained an estimate approximation in trace norm of this semigroup by the unperturbed semigroup $e^{-t\lambda'S}$. \\

In {\bf[12]}, we have regularized the operator  $H_{\mu,\lambda}$ by $\lambda''G$ where  $G = a^{*3} a^{3}$, i.e we have considered  $H_{\lambda''} = \lambda'' G + H_{\mu,\lambda}$ where $\lambda''$ is the {\it magic coupling} of Pomeron. In this case, we have established an exact relation between the degree of subordination of the non-self-adjoint perturbation operator $H_{\mu,\lambda}$ to the unperturbed operator $G$ and the number of corrections necessary for the existence of finite formula of the regularized trace.\\

The goal of this article consists to study the trace of the semigroup $e^{-tH_{\lambda''}}$, in particular to give an asymptotic expansion of this trace as $t \rightarrow 0^{+} $.}}{\color{blue}{\hf}}\\

\begin{center}

{\color{blue}{\bf 1. Introduction}}

\end{center}

Usually, quantum Hamiltonians are constructed as self-adjoint operators; for certain situations,
however, non-self-adjoint Hamiltonians are also of importance. In particular, the reggeon field
theory (as invented (1967) by V. Gribov {\bf[7]}) for the high energy behaviour of soft processes is governed by the non-self-adjoint Gribov operator\\

$H_{\lambda'',\lambda',\mu,\lambda} = \lambda''a^{*3}a^{3} + \lambda'a^{*2}a^{2} + \mu a^{*}a +i\lambda a^{*}( a + a^{*})a $\\

\n where $a$ and its adjoint $a^{*}$ are annihilation and creation operators, respectively, satisfying the canonical commutation relations $[a, a^{*}] = I$.\\

It is convenient to regard the above operators as acting on Bargmann space $E$ {\bf [4]} :\\

$ E = \{\phi: I\!\!\!\!C \rightarrow I\!\!\!\!C entire ;  \displaystyle{\int_{I\!\!\!\!C}}\displaystyle{\mid \phi(z)\mid^{2}e^{-\mid z \mid ^{2}} dxdy} < \infty \}$ \\

The Bargmann space $E$ with the paring :\\

$ < \phi, \psi > = \displaystyle{\int_{I\!\!\!\!C}}\displaystyle{\phi(z)\bar{\psi(z)}e^{-\mid z \mid ^{2}} dxdy}$\\

is a Hilbert space and $e_{n}(z) = \frac{z^n}{\sqrt{n!}}; n = 0, 1, ....$ is an orthonormal basis in $E$.\\

In this representation, the standard Bose annihilation and creation operators are defined by\\

\n $\left\{
  \begin{array}{ c }
a\phi(z) = \phi^{'}(z)\quad \quad \quad \quad \quad \quad \quad \quad \quad\quad \quad \quad \quad \quad \quad \quad \quad \quad \quad \quad \quad \quad \quad \quad\\
with \quad maximal \quad domain \quad \quad \quad \quad \quad \quad \quad \quad \quad\quad \quad \quad \quad \quad \quad \quad \quad\\
D(a) = \{\phi \in E ; \quad a\phi \in E\}\quad \quad \quad \quad \quad \quad \quad \quad \quad \quad \quad \quad \quad \quad \quad \quad \quad \quad \quad\\
\end{array}\right.$
\\
\quad\\
and\\
\quad\\
$\left\{
  \begin{array}{ c }
a^{*}\phi(z) = z\phi(z)\quad \quad \quad \quad \quad \quad \quad \quad \quad\quad \quad \quad\quad \quad \quad \quad \quad \quad \quad \quad \quad \quad \\
with \quad maximal \quad domain \quad \quad \quad \quad \quad \quad \quad \quad \quad\quad \quad \quad \quad \quad \quad \quad \quad\\
D(a^{*}) = \{\phi \in E ; \quad a^{*}\phi \in E\}\quad \quad \quad \quad \quad \quad \quad \quad \quad \quad \quad \quad \quad \quad \quad \quad\\
\end{array}\right.$

\quad\\

Notice that $D(a) = D(a^{*})$ and $D(a)\hookrightarrow E$ is compact but $\rho(a) = I\!\!\!\!C$.\\

It has been established in {\bf[8]} that $ T_{1} = aa^{*^{2}}$ is a chaotic operator and it has been shown in {\bf[9]} that $ T_{1} + T_{1}^{*}$ is symmetric but not selfadjoint.\\

Notice also that for $H_{\mu,\lambda} =  \mu a^{*}a + i\lambda a^{*}( a + a^{*})a $ with domain $\{\phi \in E$ such $\quad$ that $H_{\mu,\lambda}\phi \in E\}$, we have:\\

i) For $\mu \neq 0$ and $\lambda \neq 0$,  $H_{\mu,\lambda}$ is very far from normal and not only its self-adjoint and skew-adjoint parts do not commute but there is no inclusion in either way between their domains or with the domain of their commutator (see {\bf[13]}).\\

ii) For $\mu > 0$ and $\lambda \in I\!\!R$, $e^{-tH_{\mu,\lambda}}$ is compact (see {\bf[10]}).\\

iii) For $\mu > 0$ and $\lambda \in I\!\!R$, it has been established in {\bf [1]} that the resolvent of $H_{\mu,\lambda}$ belongs to the class $C_{1 + \epsilon} \quad \forall$ $\epsilon > 0$.\\
We recalling that a compact operator $K$ acting on a  complex Hilbert space $\mathbb{E}$ belongs to the Carlemen class Carleman $C_{p}$ of order $p$ if $\displaystyle{\sum_{n=1}^{\infty} s_{n}^{p} < \infty}$, where $s_{n}$ are s-numbers of operator $K$ i.e,the eigenvalues of the operator $\sqrt{K^{*}K}$.\\
In particular, the operator $K$ is called nuclear operator if $K \in C_{1}$ and Hilbert-Shmidt operator if $K \in C_{2}$\\
For $ p \geq 1$ the value $\displaystyle{(\sum_{n=1}^{\infty} s_{n}^{p})^{\frac{1}{p}}}$ is a norm denoted by $\mid\mid . \mid\mid_{p} $ and for $ p = 1$ it is called nuclear norm or trace norm.\\

We can consult {\bf[6]} for a systematic study of operators of Carleman class  $C_{p}$ of order $p$ .\\

Let $ H_{\lambda'} = \lambda' S +  H_{\mu,\lambda}$ where $S = a^{*^{2}}a^{2}$, this operator is more regular that $H_{\mu,\lambda}$, its semigroup $e^{-tH_{\lambda'}}$ is analytic and it has been established in {\bf[11]} that the convergence of the usual Trotter product formula for $ H_{\lambda'}$ is of classical type and can be lifted to trace-norm convergence.\\

Morever, there exist $t_{0} > 0 $ and $C > 0$ such that :\\

$\mid\mid (e^{-\frac{t}{n}\lambda'S} e^{-\frac{t}{n}H_{\mu,\lambda}})^{n} - e^{-tH_{\lambda'}}\mid\mid_{1} \leq$ $C\frac{log n}{n}$ , $ n = 2, 3, ..... $ $\forall$ $t \geq t_{0}$.\\

Now, if $H_{\mu,\lambda}$ is regularized by $\lambda''G$ where $G = a^{*^{3}}a^{3}$ and $\lambda''$ is the magic coupling of Pomeron, we can consider:\\

$ H_{\lambda^{''}} = \lambda''a^{*^{3}}a^{3} + \mu a^{*}a + i\lambda a^{*}(a  + a^{*})a$ $\hfill { } (1.1)$\\

$ \left\{\begin{array}{ c }
H_{\lambda^{''}}\phi(z)= \lambda^{''} z^{3}\phi^{'''}(z)+i\lambda z\phi^{''}(z) + (i\lambda z^{2} + \mu z)\phi^{'}(z)\quad \quad \quad \quad \quad \quad \quad \quad \quad \quad \\
with \quad maximal \quad domain \quad \quad \quad \quad \quad \quad \quad \quad \quad\quad \quad \quad \quad \quad \quad \quad \quad \quad \quad \quad \quad\\
D(H_{\lambda^{''}}) = \{\phi \in E; \quad H_{\lambda^{''}}\phi \in E\}\quad \quad \quad \quad \quad \quad \quad \quad \quad \quad \quad \quad \quad \quad \quad \quad \quad \quad \quad \quad\\
\end{array}\right.$
\quad\\

 In {\bf [2]} (see theorem 3.3, p. 595), Aimar et al  have shown that the spectrum of $H_{\lambda''}$ is discrete and that the system of generalized eigenvectors of this operator is an unconditional basis in Bargmann space $E$.\\

Recently in {\bf [12]}, we have established a regularized trace formula for $H_{\lambda''}$.\\

More precisely, we have shown the following result:\\

{\color{red}{{\bf Theorem 1.1 (see [12])}}}  \\

{\it Let $ E $ be the Bargmann space and $ H_{\lambda^{''}} = \lambda^{''}G + H_{\mu,\lambda}$ acting on $E$\\

 where\\

- $ G = a^{*3}a^{3}$ and $ H_{\mu,\lambda} = \mu a^{*}a + i\lambda a^{*}( a + a^{*})a$\\

- $a$ and $a^{*}$ are the standard Bose annihilation and creation operators satisfying \\

the commutation relation $[a , a^{*}] = I$.\\

Then \\

There exists an increasing sequence of radius $r_{n}$ such that $r_{n} \rightarrow \infty$ as $ n \rightarrow \infty$ \\

and \\

 $\displaystyle{Lim \sum_{k=0}^{n}(\sigma_{k} - \lambda''\lambda_{k}) =}$\\

 $ \displaystyle{- Lim\frac{1}{2i\pi}\int_{\gamma_{n}} Tr[\sum_{j=1}^{4}\frac{(-1)^{j-1}}{j}[H_{\mu,\lambda}(\lambda{''}G -\sigma I)^{-1}]^{j}]d\sigma}$;$ n \rightarrow \infty$ $\hfill {(1.2) } $\\

where\\

 - $\sigma_{k}$ are the eigenvalues of the operator $ H_{\lambda^{''}} = \lambda{''}G + H_{\mu,\lambda}$\\

 - $\lambda_{k} = k(k-1)(k-2)$ are the eigenvalues of the operator of $ G $ associated to eigenvectors $e_{k}(z)$\\

 - $ (\lambda{''}G - \sigma I)^{-1}$  is the resolvent of the operator $\lambda{''}G$ \\

 and\\

 - $\gamma_{n}$ is the circle of radius $r_{n}$ centered at zero in complex plane.}{\color{blue}{\hf}}\\

The goal of this article consists to study the trace of the semigroup $e^{-tH_{\lambda''}}$, in particular, to give an asymptotic expansion of this trace as $t \rightarrow 0^{+} $.\\

Our procedure consists to prove that\\

 i) The semigroup $e^{-tG}$ generated by the operator $G$ is analytic and nuclear (Gibbs analytic semigroup).\\

 ii) $ \forall \epsilon > 0, \exists C_{\epsilon} > 0  ; \mid\mid H_{\mu,\lambda}\phi\mid\mid \leq \epsilon\mid\mid G\phi\mid\mid + C_{\epsilon}\mid\mid\phi\mid\mid \forall \phi \in D(G).$ $\hfill {(1.3) } $\\

 where $D(G) = \{\phi \in E, G\phi \in E\}$.\\

 iii) The operator $H_{\mu,\lambda}G^{-\delta}$ is bounded \quad $\forall$ $\delta \geq \frac{1}{2}$.\\

 Now, with the aid of the results of Angelscu et al {\bf[3]} or of Zagrebnov {\bf [17-18 ]} with Ginibre-Gruber inequality {\bf[5]}, it easy to prove that the series of general term $S_{k}(t)_{k\in I\!\!N}$ defined by:\\

 $S_{0}(t)\phi = e^{-t\lambda''G}\phi$\\

and\\

$S_{k+1}(t)\phi = - \displaystyle{ \int_{0}^{t}e^{-(t-s)\lambda''G} H_{\mu,\lambda} S_{k}(t)\phi ds}$\\

converges to $e^{-tH_{\lambda''}}$ (nuclear norm).\\

By using the properties i), ii) and iii) we get:\\

$\mid\mid e^{-tH_{\lambda''}} - e^{-t\lambda''G} \mid\mid_{1} = t\mid\mid e^{-t\lambda''G} H_{\mu,\lambda}\mid\mid_{1} + \mid\mid (\lambda''G)^{\delta} e^{-\frac{t}{3}\lambda''G }\mid\mid_{1}O(t^{2})$; $t \rightarrow 0^{+} $ $\hfill {(1.4) } $\\

\n To establish the above results, we give in section 2 some spectral properties of semigroups $e^{-t\lambda''G}$ and $e^{-t(\lambda''G + H_{\mu,\lambda})} $ in $C_{p}$. In section 3, we give the proof of the above formula (1.4) for the trace of $e^{-tH_{\lambda''}}$ as $t \rightarrow 0^{+} $.\\
\newpage
\begin{center}
\color{blue}{\bf{ 2. Some spectral properties of semigroups $e^{-t\lambda''G}$ and $e^{-t(\lambda''G + H_{\mu,\lambda})}$ \quad in \quad Carleman \quad class \quad $C_{p}$ }}\\
\end{center}

\quad

We begin this section by given some elementary spectral properties of $G$ and $H_{\lambda{''}}$ :\\

{\color{red}{{\bf Lemma 2.1}}}\\

{\it 1) The operator $G$ is self adjoint with compact resolvent.\\

2) The eigenvalues of $G$ are $\lambda_{n} = n(n-1)(n-2)$ for $ n \geq 0 $ associated to \\

eigenvectors $e_{n}(z) = \frac{z^{n}}{\sqrt{n!}}$.\\

3) $Lim \mid\mid Ge^{-tG}\mid\mid = \frac{1}{e}$ as $ t \rightarrow 0.$\\

4) the resolvent of $G$ belongs to $ C_{p}$ \quad $\forall p > \frac{1}{3}$.\\

5) $ e^{-tG}$ is nuclear semigroup and $\mid\mid e^{-tG}\mid\mid_{1} \leq Ct^{-\frac{1}{3}}.$\\

where the constant $C$ does not depend on $t$\\

6) $ e^{-tG} \in C_{p} $ \quad $\forall p > 0$.}{\color{blue}{\hf}}\\

{\color{red}{\bf proof}}\\

1) It follows easily from Rellich theorem : \\

{\color{red}{{\bf Theorem 2.2.(see {\bf[15]}, p. 386)}}}\\

{\it Let $B$ be a self adjoint operator in a complex Hilbert space $\mathbb{E}$ such that \\ $< B\phi, \phi> \geq <\phi, \phi>, \phi \in D(B)$, where $D(B) = \{\phi \in \mathbb{E}; B\phi \in \mathbb{E}\}$. \\
Then $B$ is discrete if and only if $\{\phi \in D(B)$; $ <B\phi, \phi) \leq 1 \}$ is pre-compact.}\\

that the operator $G$ is discrete.\\

\n Also, we can prove the above result by using the following observation: \\

Since $D(a)$ is compactly embedded in Bargmann space $E$ and $D(a^{*3}a^{3}) \hookrightarrow D(a)$ is continuous then $D(a^{*3}a^{3})$ is compactly embedded in Bargmann space $E$.\\

And since $ (a^{*3}a^{3} + I) $ is invertible  then the operator $G$ is discrete .\\

2) It is evident.\\

3) As $Ge^{-tG}e_{n} = n(n-1)(n-2)e^{-tn(n-1)(n-2)}e_{n}$ we get\\

 $\mid\mid Ge^{-tG}\mid\mid  = \frac{1}{t}\quad_{_{n}}\!\!\!\!\!Sup\quad tn(n-1)(n-2)e^{-tn(n-1)(n-2)} = \frac{1}{e}$ this implies \\

$Lim \mid\mid Ge^{-tG}\mid\mid = \frac{1}{e}$ as $ t \rightarrow 0.$ $\hfill {(2.1) } $\\

4) As $G$ is self adjoint and its eigenvalues are :\\

$\lambda_{n} = n(n-1)(n-2) \sim n^{3}$ $\hfill {(2.2) } $\\

\n this implies that the series of general term  $\frac{1}{n^{3p}}$ converges $\forall p > \frac{1}{3}$ and consequently, the resolvent of  $G$ belongs to Carleman class $ C_{p}$ $\forall p > \frac{1}{3}$.\\

5) $e^{-tG}$ is self adjoint and it is of trace class, because the series of general term $e^{-t\lambda_{n}}$ converges $\forall t > 0 $.\\

Now let $ x \in [0, +\infty[$ and if $t \in [0, +\infty[ $, consider the function $f(x) = e^{-tx^{3}}$ and its derivative $f'(x) = -3tx^{2}e^{-tx^{3}}$  which non positive then the function $f(x)$ is decreasing and we have \\

$\displaystyle{\sum_{n=1}^{\infty}e^{-tn^{3}} \leq \int_{0}^{\infty}e^{-tx^{3}}dx}$ \\

By a change of variable in the above integral, we obtain that\\

$\displaystyle{\int_{0}^{\infty}e^{-tx^{3}}dx = Ct^{\frac{-1}{3}}}$\\

where the constant $C$ does not depend on $t$\\

This implies that\\

$\mid\mid e^{-tG} \mid\mid_{1} \leq \displaystyle{\sum_{n=1}^{\infty}e^{-t^{3}n} \leq  Ct^{\frac{-1}{3}}}$\\

6) $e^{-tG}$ is Carleman class $C_{p}$ $\forall p > 0 $, because the series of general term $e^{-tp\lambda_{n}}$  converges $\forall t > 0$ and $ \forall p > 0$.{\color{blue}{\hf}}\\
\newpage
{\color{red}{\bf Remark 2.3}}\\

{\it We can derive the property 6) from the fact that $e^{-tG} \in C_{1}$. In fact, since $C_{1}\subset C_{p}$ then $e^{-tG} \in C_{p}$ $\forall p > 1$.\\

Now if $ p < 1$, we choose an integer $ n$ such that $\frac{1}{n} < p$ , then  $C_{\frac{1}{n}} \subset C_{p}$.\\

Since $e^{-\frac{t}{n}G} \in C_{1}$ $\forall$ $\tau = \frac{t}{n} > 0$ then $(e^{-\frac{t}{n}G})^{n} \in C_{\frac{1}{n}}$ \\

and we obtain that\\

$ e^{-tG} \in C_{p} $ $\forall$ $ p > 0$.}{\color{blue}{\hf}}\\

From the above remark, we are now ready to prove following lemma :\\

{\color{red}{{\bf Lemma 2.4.}}}\\

{\it 1) Let $T(t)$ be a semigroup on a complex Hilbert space $\mathbb{E}$. We assume that $T(t)$ is {\color{red}{selfadjoint}} and there exists $p_{0} > 0$ such that $T(t) \in C_{p_{0}}$\quad $\forall$ $ t > 0$.\\

 Then $T(t) \in C_{p}$ \quad $\forall$ $ t > 0$, \quad $\forall$ $ p > 0$.\\

 2) Let $T(t)$ be a semigroup on a complex Hilbert space $\mathbb{E}$. We assume that there exists $p_{0} > 0$ such that  $T(t) \in C_{p_{0}}$ \quad $\forall$ $ t > 0$.\\

  Then $T(t) \in C_{p}$ \quad $\forall$ $ t > 0$, \quad $\forall$ $ p > 0$.}{\color{blue}{\hf}}\\

 {\color{red}{\bf Proof}}\\

 1) let $ t = \frac{\tau}{p_{0}} $ and  $\hat{T}(\tau) = T(\frac{\tau}{p_{0}})$ then $T(t) \in C_{p_{0}}$ if and only if $\hat{T}(\tau) \in  C_{1}$. Since  $C_{1}\subset C_{p}$ for all $ p > 1$ it follows that $\hat{T}(\tau) \in C_{p}$ for all $ p > 1$.\\

 Let $ 0 < p < 1$, Since  $\hat{T}(\tau)$ is self adjoint, then for each fixed  $\tau > 0 $, there exist a positive
decreasing sequence  $s_{n} \in l_{1}$ and an orthonormal sequence $e_{n}$ such that $\hat{T}(\tau) = \displaystyle{\sum_{n=0}^{\infty}s_{n}e_{n}\otimes e_{n}}$.\\

 Let $p = \frac{1}{\delta}$ with $\delta >1$  and $ n$ such that $ n > \delta$, then for fixed $\frac{t}{n}$ , there exist  a positive decreasing sequence $r_{k} \in l_{1}$ and an orthonormal sequence $e_{k}$ such tat $T(\frac{t}{n}) = \displaystyle{\sum_{k=0}^{\infty}r_{k}e_{k}\otimes e_{k}}$.\\

 It follows that:\\

$T(t) = [T(\frac{t}{n})] = \displaystyle{\sum_{k=0}^{\infty}r_{k}^{n}e_{k}\otimes e_{k}}$.\\

Since, $r_{k} \in l_{1}$ then $r_{k}^{n} \in l_{\frac{1}{n}} \subset l_{p}$,i.e.  $T(t) \in C_{p}$.\\

2) If $ p > p_{0}$ then  $T(t) \in C_{p}$ for all $ t > 0$ and all $ p > p_{0}$ because $C_{p_{0}} \subset C_{p}$ for all $ p > p_{0}$.\\

If $ p < p_{0}$, we choose an integer $ n $ such that $\frac{p_{0}}{n} < p$, since $T(t) \in C_{p_{0}}$ for all $ t > 0$ then $T(\frac{t}{n}) \in C_{p_{0}}$ and consequently $ T(t) = [T(\frac{t}{n}]^{n} \in C_{\frac{p_{0}}{n}} \subset C_{p}$.{\color{blue}{\hf}}\\

{\color{red}{\bf Remark 2.5}}\\

{\it Let $T$ be a compact operator, we assume that $T$  is positif and $\mid\mid T \mid\mid \leq 1$ then there exist a positive sequence $ 0 < s_{n} < 1$ and an orthonormal sequence $e_{n}$ such that $T = \displaystyle{\sum_{n=0}^{\infty}s_{n}e_{n}\otimes e_{n}}$.\\

Let $T(t) = \displaystyle{\sum_{n=0}^{\infty}s_{n}^{t}e_{n}\otimes e_{n}}$, if we choose $s_{n} \in \cap_{p=0}^{\infty}l_{p}$ where $l_{p}$ is the space of $p$-summable sequences, then $T(t) \in C_{p}$ for all $ t > 0$ and all $ p > 0$, it follows that\\

$\mid\mid T(t) \mid\mid_{p} = \displaystyle{(\sum_{n=0}^{\infty}s_{n}^{tp})^{\frac{1}{p}}}$\\

 For $s < t $ we have $s_{n}^{t} < s_{n}^{s}$ then $\mid\mid T(t) \mid\mid_{p}  < \mid\mid T(s) \mid\mid_{p}$.\\

 By using the Beppo-Levi's theorem, we deduce that $\mid\mid T(t) \mid\mid_{p} \rightarrow \infty $ as $ t \rightarrow 0$.}\\

{\color{red}{{\bf Proposition 2.6}}}\\

 {\it Let $E_{0} = \{\phi \in E; \phi(0) = 0\}$, $\mathbb{P}_{0} = \{ p \in \mathbb{P}; p(0) = 0\}$  where $ \mathbb{P}$ is the space of polynomials\\

 and \\

 $H_{\lambda''}^{min}$ with domain $D_{min}(H_{\lambda''})$ is the closure of the restriction of $H_{\lambda''}$ on $\mathbb{P}_{0}$.\\

 Then we have:\\

{\color{red}{(a)}} $\forall$ $\epsilon > 0 $, there exists $C_{\epsilon} > 0 $ such that :\\

$\mid\mid H_{\mu,\lambda}\phi\mid\mid \leq  \epsilon \mid\mid G\phi\mid\mid + C_{\epsilon}\mid\mid \phi\mid\mid$ for all $\phi \in D(G)$.$\hfill {(2.3) } $\\

{\color{red}{(b)}} $\forall$ $\epsilon > 0 $, there exists $C_{\epsilon} > 0 $ such that :\\

$\mid < H_{\mu,\lambda}\phi, \phi >\mid \leq  \epsilon < G\phi, \phi > + C_{\epsilon}\mid\mid \phi\mid\mid^{2}$ $\forall$  $\phi \in D(G)$.$\hfill {(2.4) } $\\

{\color{red}{(c)}} For $\lambda'' > 0$ and $\forall$  $\epsilon $ ; $0 < \epsilon < \lambda'' $, there exists $C_{\epsilon} > 0 $ such that :\\

$Re < H_{\lambda''}\phi, \phi > \geq ( \lambda'' - \epsilon )< G\phi, \phi > -C_{\epsilon}\mid\mid \phi\mid\mid^{2}$ $\forall$ $\phi \in D(G)$ $\hfill {(2.5) } $\\

\n in particular, the range of $ H_{\lambda''}$ is closed.\\

{\color{red}{(d)}} For $\lambda^{''} \geq 0$ and $\mu > 0$, $H_{\lambda{''}}$ is accretive and maximal.\\

{\color{red}{(e)}} $D_{max}(H_{\lambda{''}}) = D_{min}(H_{\lambda{''}}) = D(G)$.\\

{\color{red}{(f)}} $ - H_{\lambda''}$ generates an analytic semigroup $e^{-tH_{\lambda''}}$, $ t > 0 $.\\

{\color{red}{(g)}} For $\lambda^{''} \neq 0$, the resolvent of $H_{\lambda{''}}$ belongs to Carleman class $C_{p}$ for all $ p > \frac{1}{3}$.}{\color{blue}{\hf}}\\

{\color{red}{{\bf Proof}}}\\

a) Let $\phi \in E_{0}$, we have $\phi(z) = \displaystyle{\sum_{n=1}^{\infty}a_{n}e_{n}(z)}$. Then \\

$H_{\mu,\lambda}\phi(z) = \displaystyle{\sum_{n=1}^{\infty}[\mu na_{n} + i\lambda (n-1)\sqrt{n}a_{n-1} + i\lambda n\sqrt{n+1}a_{n+1}] e_{n}(z)}$ \\

and\\

$G\phi(z) = \displaystyle{\sum_{n=1}^{\infty}n(n-1)(n-2)a_{n}e_{n}(z)}$\\

we remark that there exists $C > 0$ such that $\mid \mid H_{\mu,\lambda}\phi\mid \mid^{2} \leq C \displaystyle{\sum_{n=1}^{\infty}n^{3} \mid a_{n} \mid^{2}}$\\

and\\

$\mid \mid G\phi\mid \mid^{2} \geq \frac{1}{36} \displaystyle{\sum_{n=1}^{\infty}n^{6}\mid a_{n}\mid^{2}}$.\\

Now, by using the Young's inequality, we get:\\

$\forall \epsilon > 0$ , $k^{3} \leq \epsilon k^{6} + \frac{1}{\epsilon}$ $\forall k \in I\!\!N$.\\

\n this implies that :\\

$\forall $ $\epsilon > 0$, there exists $C_{\epsilon} > 0 $ such that $\mid\mid H_{\mu,\lambda}\phi\mid\mid \leq  \epsilon \mid\mid G\phi\mid\mid + C_{\epsilon}\mid\mid \phi\mid\mid $ for all $\phi \in D(G)$.\\

b) $ < H_{\mu,\lambda}\phi, \phi > = \mu \mid\mid a\phi\mid\mid^{2} + i\lambda< a^{2}\phi, a\phi > + i\lambda< a\phi, a^{2}\phi >$ for all $\phi \in D(H_{\mu,\lambda})$.\\

Then\\

$\mid < H_{\mu,\lambda}\phi, \phi >\mid \leq \mu \mid\mid a\phi\mid\mid^{2} + 2\mid \lambda \mid \mid \mid a\phi\mid\mid.\mid\mid a^{2}\phi \mid\mid$ pour tout $\phi \in D(H_{\mu,\lambda})$.\\

With the aid of following inequalities:\\

i) $\forall$ $\epsilon_{1} > 0$ $\mid \mid a\phi\mid\mid.\mid\mid a^{2}\phi \mid\mid \leq  \epsilon_{1}\mid\mid a^{2}\phi \mid\mid^{2} + \frac{1}{\epsilon_{1}}\mid\mid a\phi \mid\mid^{2}.$\\

ii) $\forall$ $\epsilon_{2} > 0$ there exists $C_{\epsilon_{2}} > 0 $ such that $\mid\mid a\phi \mid\mid^{2}\leq \epsilon_{2}\mid\mid a^{3}\phi \mid\mid^{2} + C_{\epsilon_{2}}\mid\mid \phi\mid\mid^{2}$\\

iii) $\mid\mid a^{2}\phi \mid\mid^{2} \leq \mid\mid a^{3}\phi \mid\mid^{2}$\\

we obtain:\\

$\forall$ $\epsilon > 0 $, there $C_{\epsilon} > 0 $ such that $\mid < H_{\mu,\lambda}\phi, \phi >\mid \leq  \epsilon < G\phi, \phi > + C_{\epsilon}\mid\mid \phi\mid\mid^{2}$ for all $\phi \in D(G)$.\\

c) Since $Re < H_{\lambda''}\phi, \phi > = \lambda''< G\phi, \phi > + < H_{\mu,\lambda}\phi, \phi >$ and $ \lambda'' > 0$\\

Then \\

$Re < H_{\lambda''}\phi, \phi > \geq \lambda''< G\phi, \phi > - \mid < H_{\mu,\lambda}\phi, \phi >\mid$.\\

By using the above property we get:\\

 $Re < H_{\lambda''}\phi, \phi > \geq \lambda''< G\phi, \phi > - \epsilon< G\phi, \phi > - C_{\epsilon}\mid\mid \phi\mid\mid^{2} = (\lambda'' -\epsilon) < G\phi, \phi > - C_{\epsilon}\mid\mid \phi\mid\mid^{2}$\\

We choose $ 0 < \epsilon < \lambda'' $ to deduce that\\

 $Re < H_{\lambda''}\phi, \phi > \geq - C_{\epsilon}\mid\mid \phi\mid\mid^{2}$\\

Consequently the range of $ H_{\lambda''}$ is closed.\\

d) Since $a^{*}(a + a^{*})a $ is symmetric operator then \\

$Re < H_{\lambda''}\phi,\phi > = \lambda''\mid \mid a^{3}\phi\mid\mid^{2} + \mu\mid\mid a\phi\mid\mid^{2} \geq \mu\mid\mid \phi\mid\mid^{2}$, $\forall \phi \in D_{min}(H_{\lambda''})$ \\

Now for $\lambda''\geq 0$ and  $\mu >0$ we deduce that :\\

$Re <H_{\lambda''}\phi,\phi > \geq \mu\mid\mid \phi\mid\mid^{2}$,$ \forall \phi \in D_{min}(H_{\lambda''})$ .\\

This inequality will be not verified if we kept constant functions in Bargmann space $E$.\\

Now, we would like to show that there exists  $\beta_{0} \in I\!\!R$ ; $H_{\lambda''} + \beta_{0}I$ is invertible.\\

We rewrite $H_{\lambda''}$ in the following form:\\

$H_{\lambda''} = \lambda''( G + \frac{1}{\lambda''}H_{\mu,\lambda}$ et $G + \beta I + \frac{1}{\lambda''}H_{\mu,\lambda} = [I + \frac{1}{\lambda''}H_{\mu,\lambda}(G + \beta I )^{-1}](G + \beta I)$.\\

Using the property a) to get:\\

$ \mid\mid \frac{1}{\lambda''}H_{\mu,\lambda}(G + \beta I)^{-1}\psi\mid\mid \leq \epsilon \mid\mid G (G + \beta I )^{-1}\psi\mid\mid + C_{\epsilon}\mid\mid (G + \beta I )^{-1}\psi\mid\mid$\\

$\leq \epsilon \mid\mid (G +\beta I - \beta I) (G + \beta I )^{-1}\psi\mid\mid + C_{\epsilon}\mid\mid (G + \beta I)^{-1}\psi\mid\mid$\\

$ \leq \epsilon \mid\mid \psi \mid\mid + (\epsilon \beta + C_{\epsilon})\mid\mid (G + \beta I)^{-1}\psi\mid\mid$\\

\n and as $\mid\mid (G + \beta I )^{-1}\mid\mid \leq \frac{1}{\beta}$ then \\

$ \mid\mid \frac{1}{\lambda''}H_{\mu,\lambda}(G + \beta I )^{-1}\psi\mid\mid \leq (2\epsilon + \frac{C_{\epsilon}}{\beta})$. \\

Now, we choose $ 0 < \epsilon < \frac{1}{2}$ and $\beta > \frac{C_{\epsilon}}{1-2\epsilon}$ to obtain :\\

$ \mid\mid \frac{1}{\lambda''}H_{\mu,\lambda}(G + \beta I )^{-1}\mid\mid < 1$.\\

this implies $H_{\lambda''}^{min} + \beta I$ is invertible.\\

e) We begin to show that  $D_{max}(H_{\lambda{''}}) = D_{min}(H_{\lambda{''}})$.\\

First,  $D_{min}(H_{\lambda{''}}) \subset  D_{max}(H_{\lambda{''}})$ is trivial.\\

To show that $D_{max}(H_{\lambda{''}}) \subset D_{min}(H_{\lambda{''}})$, $\phi \in D_{max}(H_{\lambda{''}})$ then $ (H_{\lambda{''}} + \beta I)\phi \in E_{0}$ for all $\beta$.\\

\n Since there exists  $\beta_{0} \in I\!\!\!\! R$ such that $H_{\lambda''}^{min} + \beta_{0} I$ is invertible of $D_{min}(H_{\lambda{''}})$ on $E_{0}$, then there exists
$\phi_{1} \in D_{min}(H_{\lambda{''}})$ such that :\\

$(H_{\lambda''} + \beta_{0} I)\phi = (H_{\lambda''}^{min} + \beta_{0} I)\phi_{1} $, in particular, we have $(H_{\lambda''} + \beta_{0} I)(\phi - \phi_{1}) = 0$.\\

\n To deduce that $\phi = \phi_{1}$, we need that  Ker ($H_{\lambda''} + \beta_{0} I$) = \{0\}.\\

\n we recall that the range of $H_{\lambda''}^{min} + \beta_{0} I$ is closed and the formal adjoint of $H_{\lambda''}$ is $H_{\lambda''}$ where we substitute $\lambda$ by its opposite.\\

 Now, since the adjoint of the minimal is formal adjoint of the maximal and $H_{\lambda''}^{min} + \beta_{0} I$ is invertible then $Ker (H_{\lambda''} + \beta_{0} I) = \{0\}$ this implies that $\phi = \phi_{1}$ and $\phi \in D_{min}(H_{\lambda''})$.\\

From the inequality a) and the theorem 111 in book's Kato {\bf[4]}, we deduce $D_{max}(H_{\lambda''}) = D(G)$.\\

f) From the inequality a) and the theorem 2.1 in book's Pazy {\bf[16]}, we deduce that $ - H_{\lambda''}$ generates analytic semigroup $e^{-tH_{\lambda''}}$ $ t > 0$.\\

g) For $\lambda^{''} \neq 0$ the resolvent of $H_{\lambda{''}}$ is Carleman class $C_{p}$ for all $ p > \frac{1}{3}$, this property is the lemma 4.1 of {\bf[2]}.\\ This ends the proof of this proposition.{\color{blue}{\hf}}\\

Now, the above properties d) and e) allow us to show the following theorem:\\

{\color{red}{\bf Theorem 2.7}}\\

{\it $- H_{\lambda{''}}$ generates a semigroup $e^{-tH_{\lambda''}}$ of Carleman $C_{p}$ for all $ p > 0$ and all $ t > 0$.}{\color{blue}{\hf}}\\

{\color{red}{\bf Proof}}\\

 Let $E_{0} = \{\phi \in E ; \phi(0) = 0 \}$ then on  $E_{0}$ we have:\\

 $Re<H_{\lambda''}\phi, \phi> = \lambda{''}\mid\mid A^{3}\phi\mid\mid^{2} + \mu\mid\mid A\phi\mid\mid^{2} \geq \mu \mid\mid \phi \mid\mid^{2}$ for $ \lambda{''} \geq 0$ and $\mu > 0$. \\

 From this inequality we deduce that $0 $ belongs to resolvent set $\rho(H_{\lambda''})$ of the operator $H_{\lambda''}$.\\

 Let $T(t) = \displaystyle{\int_{0}^{t}e^{-sH_{\lambda''}}}\phi ds$ then\\

 $T(t) = H_{\lambda''}^{-1}( I - e^{-tH_{\lambda''}})$ and as the resolvante of $H_{\lambda''}$ is Carleman class $C_{p}$ for all $ p > \frac{1}{3}$ and the operator $I - e^{-tH_{\lambda''}}$is bounded then $T(t)$ is Carleman $C_{p}$ for all $ p > \frac{1}{3}$ and with the aid of lemma 2.4, we end the proof.{\color{blue}{\hf}}\\

\begin{center}

{\color{blue}{\bf{ 3. Asymptotic expansion of trace of $e^{-tH_{\lambda''}}$ as $t \rightarrow 0^{+}$}}}

\end{center}

We put $S(s) = e^{-(t-s)\lambda''G} e^{-sH_{\lambda''}}$ then for $\phi \in D(G)$, the application $\phi \rightarrow S(s)\phi $ is differentiable and $S'(s)\phi = \lambda''Ge^{-(t-s)\lambda''G}   e^{-sH_{\lambda''}} - e^{-(t-s)\lambda''G}H_{\lambda''}e^{-sH_{\lambda''}} = -e^{-(t-s)\lambda''G}H_{\mu,\lambda}e^{-sH_{\lambda''}}$.\\

On $[0, t]$, we have $\displaystyle{ \int_{0}^{t}S'(s)\phi ds = S(t)\phi - S(0)\phi }= -\displaystyle{ \int_{0}^{t}e^{-(t-s)\lambda''G}H_{\mu,\lambda}e^{-sH_{\lambda''}}\phi ds}$ \\

This implies that $e^{-tH_{\lambda''}}\phi $ is solution of the following integral equation :\\

$e^{-tH_{\lambda''}}\phi - e^{-t\lambda''G}\phi = -\displaystyle{ \int_{0}^{t}e^{-(t-s)\lambda''G}H_{\mu,\lambda}e^{-sH_{\lambda''}}\phi ds}$ $\hfill {(3.1) } $\\

or\\

$e^{-tH_{\lambda''}}\phi - e^{-t\lambda''G}\phi = -\displaystyle{ \int_{0}^{t}N(t,s)e^{-sH_{\lambda''}}\phi ds}$ with $N(t,s) = e^{-(t-s)\lambda''G}H_{\mu,\lambda}$ $\hfill {(3.2) } $\\

It is well known that the solution of equation (3.1) can be obtained by successive approximation method:\\

$ e^{-tH_{\lambda''}} = \displaystyle{\sum_{k=0}^{\infty}S_{k}(t)}$ $\hfill {(3.3) } $\\

where \\

$S_{0}(t)\phi = e^{-t\lambda''G}\phi$\\

and\\

$S_{k+1}(t)\phi = - \displaystyle{ \int_{0}^{t}e^{-(t-s)\lambda''G}H_{\mu,\lambda} S_{k}(t)\phi ds}$\\

the convergence of (3.3) is in operator norm.\\

 Notice that:\\

 $S_{1}(t)\phi = - \displaystyle{ \int_{0}^{t}e^{-(t-t_{1})\lambda''G}H_{\mu,\lambda}e^{-t_{1}\lambda''G}\phi dt_{1}}$\\

\quad\\

 $S_{2}(t)\phi = \displaystyle{ \int_{0}^{t}\displaystyle{ \int_{0}^{t_{1}}e^{-(t-t_{1})\lambda''G}H_{\mu,\lambda}e^{-(t_{1}-t_{2})\lambda''G}H_{\mu,\lambda}e^{-t_{2}\lambda''G}\phi dt_{2}dt_{1}}}$\\

 \quad\\

 $S_{3}(t)\phi = (-1)^{3} \displaystyle{ \int_{0}^{t}\displaystyle{ \int_{0}^{t_{1}}}\displaystyle{ \int_{0}^{t_{2}}e^{-(t-t_{1})\lambda''G}H_{\mu,\lambda}e^{-(t_{1}-t_{2})\lambda''G}H_{\mu,\lambda}e^{-(t_{2}-t_{3})\lambda''G}H_{\mu,\lambda}e^{-t_{3}\lambda''G}\phi dt_{3}dt_{2}dt_{1}}}$\\

  \quad\\

$S_{k}(t)\phi = (-1)^{k} \displaystyle{ \int_{0}^{t}\displaystyle{ \int_{0}^{t_{1}}\displaystyle{ \int_{0}^{t_{2}} .................................}}}
\displaystyle{ \int_{0}^{t_{k-1}}}$\\

$e^{-(t-t_{1})\lambda''G}H_{\mu,\lambda}e^{-(t_{1}-t_{2})\lambda''G}H_{\mu,\lambda}e^{-(t_{2}-t_{3})\lambda''G}H_{\mu,\lambda}
............................................ $\\

$e^{-(t_{k-1}-t_{k})\lambda''G}H_{\mu,\lambda}e^{-t_{k}\lambda''G}\phi dt_{k}............................................ dt_{3}dt_{2}dt_{1}$.\\

{\color{red}{\bf Lemma 3.1}}\\

{\it The series $e^{-tH_{\lambda''}} = \displaystyle{\sum_{k = 0}^{\infty}S_{k}(t)}$ converges in trace norm.}{\color{blue}{\hf}}\\

{\color{red}{\bf Proof}}\\

The convergence in trace norm is obtained by using the results of Angelescu-Nenciu-Bundaru, in particular their  proposition in {\bf[3]}  or the results of Zagrebnov in {\bf [17-18 ]}, in particular the theorem 2.1 in {\bf [18]} with the aid of  Ginibre-Gruber's inequality in {\bf[5]}.\\

Now we are going to derive an asymptotic expansion of the trace of $e^{-tH_{\lambda''}}$ as asymptotique lorsque $t\rightarrow 0^{+}$\\

For $ s \geq 0 $ we have\\

$e^{-sH_{\lambda''}}\phi = e^{-s\lambda''G}\phi -\displaystyle{ \int_{0}^{s}N(s,s_{1})e^{-s_{1}H_{\lambda''}}\phi ds_{1}}$ with $N(s,s_{1}) = e^{-(s-s_{1})\lambda''G}H_{\mu,\lambda}$\\

and for $ 0 \leq s \leq t$, we substitute the above expression of $e^{-sH_{\lambda''}}\phi$ in (3.2) to get  :\\

$e^{-tH_{\lambda''}}\phi - e^{-t\lambda''G}\phi = -\displaystyle{ \int_{0}^{t}N(t,s)[e^{-s\lambda''G}\phi -\displaystyle{ \int_{0}^{s}N(s,s_{1})e^{-s_{1}H_{\lambda''}}\phi }]ds_{1}ds}$\\

$= -\displaystyle{\int_{0}^{t}N(t,s)e^{-s\lambda''G}\phi ds}$ + $\displaystyle{ \int_{0}^{t}\displaystyle{ \int_{0}^{s}N(t,s)N(s,s_{1})e^{-s_{1}H_{\lambda''}}\phi ds_{1}ds}}$\\

$= -\displaystyle{\int_{0}^{t}e^{-(t-s)\lambda''G}H_{\mu,\lambda}e^{-s\lambda''G}\phi ds}$ + $\displaystyle{\int_{0}^{t}\displaystyle{ \int_{0}^{s}e^{-(t-s)\lambda''G}H_{\mu,\lambda}e^{-(s-s_{1})\lambda''G}H_{\mu,\lambda}e^{-s_{1}H_{\lambda''}}\phi ds_{1}ds}}$\\

Then we have: \\

$e^{-tH_{\lambda''}}\phi - e^{-t\lambda''G}\phi = \displaystyle {\int_{0}^{t}e^{-(t-s)\lambda''G}H_{\mu,\lambda}e^{-s\lambda''G}\phi ds}$ + $\displaystyle {\int_{\Delta}e^{-(t-s)\lambda''G}H_{\mu,\lambda}e^{-(s-s_{1})\lambda''G}H_{\mu,\lambda}e^{-s_{1}H_{\lambda''}}\phi ds_{1}ds}$ $ \hfill {(3.3) } $\\

where $\Delta$ is the triangle $\{(s_{1},s); 0 \leq s_{1} \leq s \leq t \}$.\\

The matrix associated to $H_{\mu,\lambda}$ in the basis $ \{e_{n}\}$ can be written in this form :\\

$H_{\mu,\lambda} e_{n} = i\lambda(n-1)\sqrt{n}e_{n-1} + n\mu e_{n} + i\lambda n\sqrt{n+1}e_{n+1}$  $\hfill {(3.4) } $\\

- The  family of infinite matrices associated to $H_{\mu,\lambda}$ is tridiagonal of the form $J + i\lambda H^{(I)}$, where the matrix $J$ is diagonal with entries $J_{nn} := n\mu$ , and the matrix $H^{(i)}$ is off-diagonal, with nonzero entries $H_{n,n+1}^{(i)} = H_{n+1,n}^{(i)} := H_{n}^{(I)} = n\sqrt{n+1}$.\\

- The  family of infinite matrices associated to $\lambda''G$ is diagonal with entries $G_{nn} := \lambda''n(n-1)(n-2)$\\

Let the infinite matrix $ ^{t\!}H_{\mu,\lambda} $ be obtained from $H_{\mu,\lambda}$ by transposing of the elements and the infinite matrix $ H_{\mu,\lambda}^{\bot}$ be obtained from $H_{\mu,\lambda}$ by transposing and by taking complex conjugates of the elements. Then observe that\\

i) $H_{\mu,\lambda}$ is symmetric complex matrix i.e. $ H_{\mu,\lambda} =\quad ^{t\!}H_{\mu,\lambda}$.\\

ii) $ H_{\mu,\lambda} \neq  H_{\mu,\lambda}^{\bot}$ (The symbol $\bot$ represents Dirac Hermitian conjugation;that is, transpose and complex conjugate.)\\

iii) As $H_{n}^{(I)} = O(n^{\alpha})$ with $\alpha = \frac{3}{2} > 1$ then the standard perturbation theory is not applicable.\\

iv) As $G_{nn} = \lambda''n(n-1)(n-2)$ then $G_{nn} = O(n^{3})$ as $n\rightarrow \infty$\\

v) For other properties on the matrix associated to $H_{\mu,\lambda}$, we can consult {\bf [13]}.\\

vi) from the above observations, we deduce that the operator $H_{\mu,\lambda}G^{-\delta}$ is bounded for all $ \delta \geq \frac{1}{2}$.\\

Now, we present the aim result of this work in following theorem:\\

{\color{red}{\bf Theorem 3.2 }}\\

{\it Let  $H_{\lambda''} = \lambda'' G + H_{\mu,\lambda}$ the Gribov's operator acting on Bargmann's space.\\

where\\

$ G = a^{*3} a^{3}$ and $ H_{\mu,\lambda} =\mu a^* a + i \lambda a^* (a+a^*)a$\\

$[a, a^{*}] = I$ and $(\lambda'', \mu,\lambda)$ are reel parameters and $i^{2} = -1$.\\

Then\\

$\mid\mid e^{-tH_{\lambda''}} - e^{-t\lambda''G} \mid\mid_{1} = t\mid\mid e^{-t\lambda''G} H_{\mu,\lambda}\mid\mid_{1} + \mid\mid (\lambda''G)^{\delta} e^{-\frac{t}{3}\lambda''G }\mid\mid_{1}O(t^{2})$.}{\color{blue}{\hf}}\\

{\color{red}{\bf Proof }}\\

{\color{red}{\bf a}}) We begin by computing the trace of the operator:\\

$I_{1}(t) = \displaystyle {\int_{0}^{t}e^{-(t-s)\lambda''G}H_{\mu,\lambda}e^{-s\lambda''G}ds}$\\

We have\\

$\mid\mid I_{1}(t)\mid\mid_{1} = \displaystyle {\int_{0}^{t}\mid\mid e^{-(t-s)\lambda''G}H_{\mu,\lambda}e^{-s\lambda''G}\mid\mid_{1}ds}$\\

$= \displaystyle {\int_{0}^{t}\displaystyle {\sum_{n=1}^{\infty}<e^{-(t-s)\lambda''G}H_{\mu,\lambda}e^{-s\lambda''G}e_{n},e_{n}> ds}}$\\

$ = \displaystyle {\int_{0}^{t}\displaystyle{\sum_{n=1}^{\infty} < e^{-t\lambda''G}H_{\mu,\lambda}e^{-s\lambda_{n}}e_{n}, e^{s\lambda_{n}}e_{n} > ds}}$ because $e^{s\lambda''G}$ is self adjoint\\

$= \displaystyle {\int_{0}^{t}\displaystyle {\sum_{n=1}^{\infty}<e^{-t\lambda''G}H_{\mu,\lambda}e_{n}, e_{n}>ds}}$\\

$ = \displaystyle {\int_{0}^{t}\mid\mid e^{-t\lambda''G}H_{\mu,\lambda}\mid\mid_{1}ds}$\\

$ = t\mid\mid e^{-t\lambda''G}H_{\mu,\lambda}\mid\mid_{1}$\\

Then we deduce that:\\

$\mid\mid I_{1}(t)\mid\mid_{1} = t\mid\mid e^{-t\lambda''G}H_{\mu,\lambda}\mid\mid_{1}$ $\hfill {(3.5) } $\\

{\color{red}{\bf b}}) We begin to recall the symmetry property of the norm in Carleman class $C_{p}$\\

The symmetry of the norm in $C_{p}$ means that\\

$\mid\mid K_{1}K_{2}K_{3}\mid\mid_{p} \leq \mid\mid K_{1}\mid\mid.\mid\mid K_{2}\mid\mid_{p}\mid\mid K_{3}\mid\mid.$ $\hfill {(3.6) } $\\

for any bounded operators $K_{1}$ and $K_{3}$ and $K_{2}\in C_{p}$\\

Consider the trace of the operator \\

$I_{2}(t) = \displaystyle { \int_{\Delta}e^{-(t-s)\lambda''G}H_{\mu,\lambda}e^{-(s-s_{1})\lambda''G}H_{\mu,\lambda}e^{-s_{1}H_{\lambda''}}ds_{1}ds}$\\

and let $\delta \geq\frac{1}{2}$ such that $H_{\mu,\lambda}G^{-\delta}$ bounded, then \\

$\mid\mid I_{2}(t)\mid\mid_{1} = $

$\displaystyle { \int_{\Delta}\mid\mid e^{-(t-s)\lambda''G}H_{\mu,\lambda}G^{-\delta}[G^{\delta}e^{-(s-s_{1})\lambda''G}]H_{\mu,\lambda}G^{-\delta}[G^{\delta}e^{-s_{1}H_{\lambda''}}]\mid\mid_{1} ds_{1}ds}$\\

As $t$ can be written as sum of three positif numbers $t = (t-s) +(s-s_{1}) + s_{1}$. It follows that at least one of them is not less than $\frac{t}{3}$; suppose, for example, that $ s-s_{1} \geq \frac{t}{3}$\\

Then \\

$\mid\mid G^{\delta}e^{-(s-s_{1})\lambda''G}\mid\mid_{1} \leq$ $\mid\mid G^{\delta}e^{-\frac{t}{3}\lambda''G}\mid\mid_{1}$\\

By using the inequality (3.6) we deduce that\\

$\mid\mid e^{-(t-s)\lambda''G}H_{\mu,\lambda}G^{-\delta}[G^{\delta}e^{-(s-s_{1})\lambda''G}]H_{\mu,\lambda}G^{-\delta}[G^{\delta}e^{-s_{1}H_{\lambda''}}]\mid\mid_{1}$
$ \leq \mid\mid e^{-(t-s)\lambda''G}H_{\mu,\lambda}G^{-\delta}\mid\mid . \mid\mid G^{\delta}e^{-(s-s_{1})\lambda''G}\mid\mid_{1}.
 \mid\mid H_{\mu,\lambda}G^{-\delta}[G^{\delta}e^{-s_{1}H_{\lambda''}}\mid\mid$\\

$ \leq \mid\mid H_{\mu,\lambda}G^{-\delta} \mid\mid^{2}\mid\mid G^{\delta}e^{-\frac{t}{3}\lambda''G}\mid\mid_{1}$\\

Then we have\\

 $\mid\mid I_{2}(t)\mid\mid_{1} \leq \mid\mid H_{\mu,\lambda}G^{-\delta} \mid\mid^{2}\mid\mid G^{\delta}e^{-\frac{t}{3}\lambda''G}\mid\mid_{1}\displaystyle { \int_{\Delta}dsds_{1}}$\\

 $\leq \mid\mid H_{\mu,\lambda}G^{-\delta} \mid\mid^{2}\mid\mid G^{\delta}e^{-\frac{t}{3}\lambda''G}\mid\mid_{1}t^{2}$\\

It follows that\\

$\mid\mid I_{2}(t)\mid\mid_{1} = \mid\mid G^{\delta}e^{-\frac{t}{3}\lambda''G}\mid\mid_{1}O(t^{2})$ $\hfill {(3.7) } $\\

and consequently we have\\

$\mid\mid e^{-tH_{\lambda''}} - e^{-t\lambda''G} \mid\mid_{1} = t\mid\mid e^{-t\lambda''G} H_{\mu,\lambda}\mid\mid_{1} + \mid\mid (\lambda''G)^{\delta} e^{-\frac{t}{3}\lambda''G} \mid\mid_{1}$  $O(t^{2})$ $\hfill {(3.8)}${\color{blue}{\hf}}\\

\begin{center}
{\large {\bf References}}\\
\end{center}

\n {\bf[1]} M.T. Aimar, A. Intissar, J.-M. Paoli, Quelques nouvelles propri\'et\'es de r\'egularit\'e de l'op\'erateur de Gribov, Comm. Math. Phys. 172 (1995) 461-466.\\

\n {\bf [2]} M.T. Aimar, A. Intissar, A. Jeribi, On an unconditional basis of generalized eigenvectors of the nonself-adjoint Gribov Operator in Bargmann Space,Journal of Mathematical Analysis and Applications 231, (1999), 588-602.\\

\n {\bf[3]} N. Angelescu, G. Nenciu, M. Bundaru, On the perturbation of Gibbs semigroups, Comm. Math. Phys., 42 (1975).\\

\n {\bf[4]} V. Bargmann, On a Hilbert space of analytic functions and an associated integral transform I, Comm. Pure
Appl. Math. 14 (1962) 187-214.\\

\n {\bf[5]} J. Ginibre, C. Gruber, Green functions of anisotropic Heisenberg model, Comm. Math. Phys. 11 (1969).\\

\n{\bf[6]} I.C. Gohberg, M. Krein, Introduction to the theory of linear non-self adjoint operators, 18,Providence R.I: A.M.S., (1969).\\

\n {\bf[7]} V. Gribov, A reggeon diagram technique, Soviet Phys. JETP 26 (1968), no. 2, 414-423.\\

\n {\bf[8]} A. Intissar, On a chaotic weighted Shift $z^{p}d^{p+1}/dz^{p+1}$ of order p in Bargmann space, Advances in Mathematical Physics, Article ID 471314, (2011).\\

\n{\bf[9]} A. Intissar, Analyse de Scattering d'un op\'erateur cubique de Heun dans l'espace de Bargmann, Comm. Math. Phys., 199 (1998) 243-256.\\

\n {\bf[10]} A. Intissar, Etude spectrale d'une famille d'op\'erateurs non-sym\'etriques intervenant dans la th\'eorie des champs de Reggeons, Comm. Math. Phys. 113 (1987) 263-297.\\

\n{\bf [11]} A. Intissar, Approximation of the semigroup generated by the Hamiltonian of Reggeon field theory in Bargmann space, Journal of Mathematical Analysis and Applications, vol. 305, no. 2,(2005), pp. 669-689\\

\n{\bf[12]} A. Intissar, Regularized trace of magic Gribov operator on Bargmann space, Journal of Mathematical Analysis and Applications (submitted), {\color{blue} arXiv:1311.1394} \\

\n {\bf [13]} A. Intissar, Analyse Fonctionnelle et Th$\acute{e}$orie Spectrale pour les Op$\acute{e}$rateurs Compacts Non Auto-Adjoints, Editions Cepadues, Toulouse, (1997).\\

\n{\bf[14]} T. Kato, Perturbation Theory for Linear Operators, Springer, Berlin, (1966).\\

\n {\bf [15]} M.A. Naymark,  Linear differential operators. Nauka, M. 528 (1969)\\
\n
{\bf [16]} A. Pazy,  Semigroups of linear operators and applications to partial differential equations, Springer-Verlag New York, Inc. (1983)\\

\n{\bf[17]} V.A. Zagrebnov, On the families of Gibbs semigroups, Commun. Math. Phys. 76 (1980) 269-276\\

\n {\bf[18]} V.A. Zagrebnov, Perturbations of Gibbs semigroups, Commun. Math. Phys. 120 (1989) 653-664\\

\end{document}